\documentclass[11pt]{article}

\input{amssym.def}
\input{amssym}
\usepackage{eufrak}

\usepackage{graphicx}

\newtheorem{theorem}{Theorem}

\newtheorem{definition}[theorem]{Definition}

\newtheorem{proposition}[theorem]{Proposition}
\newtheorem{remark}[theorem]{Remark}

\def\PP{\cal P}

\def\id{\rm{id}}
\def\tx{\tilde{x}}
\def\ggg{{\frak{g}}}

\def\qq{q^{-1}}

\def\h{\hbar}

\def\Tr{\rm Tr}

\def\LL{{\cal{L}}}

\def\de{\delta}

\def\ot{\otimes}

\def\C{{\Bbb C}}

\def\Z{{\Bbb Z}}

\def\End{{\rm End}}
\def\vv{V^{\otimes 2}}

\def\RR{{\cal R}}

\def\RRR{{\bf{R}}}

\def\al{{\alpha}}

\def\be{\begin{equation}}
\def\ee{\end{equation}}

\def\bea{\begin{eqnarray}}
\def\eea{\end{eqnarray}}

\begin{document}

\makeatletter
\renewcommand{\theequation}{{\thesection}.{\arabic{equation}}}
\@addtoreset{equation}{section} \makeatother

\title{Quantum doubles of Fock type and bosonization}
\author{\rule{0pt}{7mm} Dimitry Gurevich\thanks{gurevich@ihes.fr}\\
{\small\it Institute for Information Transmission Problems}\\
{\small\it Bolshoy Karetny per. 19,  Moscow 127051, Russian Federation}\\
\rule{0pt}{7mm} Pavel Saponov\thanks{Pavel.Saponov@ihep.ru}\\
{\small\it
National Research University Higher School of Economics,}\\
{\small\it 20 Myasnitskaya Ulitsa, Moscow 101000, Russian Federation}\\
{\small \it and}\\
{\small \it
Institute for High Energy Physics, NRC "Kurchatov Institute"}\\
{\small \it Protvino 142281, Russian Federation}}

\maketitle

\begin{abstract}
We introduce analogs of  creation and annihilation operators, related to involutive and Hecke symmetries $R$,  and perform bosonic and fermionic realization  of the
modified Reflection Equation algebras in terms of  the so-called  Quantum Doubles of Fock type. Also, we introduce Quantum Doubles of Fock type, associated with
Birman-Murakami-Wenzl symmetries  coming from orthogonal or simplectic Quantum Groups  and exhibit the  algebras obtained by means of the corresponding
bosonization (fermionization). Besides, we apply  this scheme to current braidings arising from Hecke symmetries $R$ via the Baxterization procedure.
\end{abstract}

{\bf AMS Mathematics Subject Classification, 2010:} 81R50


\section{Introduction}
By an associative double  we mean a couple $(A,B)$ of associative algebras $A$ and $B$ endowed with the so-called permutation map
\be
\sigma :A \ot B\to B\ot A,
\label{PM}
\ee
which satisfies certain requirements (see Section 3).

By a quantum double (QD) we mean an associative double where the map $\sigma$ is defined by means of a braiding --- constant or current (i.e. depending on spectral
parameters), different from a (super-)flip. A typical example is the so-called Heisenberg double constructed of an RTT algebra and a reflection equation (RE) one, associated
with the same braiding $R$. A number of other examples of QD can be found in \cite{GPS2}.

By a bosonic QD  of Fock type we mean a double composed of the algebras $\mathrm{Sym}_R(V)$ and $\mathrm{Sym}_R(V^*)$, where $R:\vv\to \vv$ is a braiding, while
 $\mathrm{Sym}_R(V)$ and $\mathrm{Sym}_R(V^*)$ are $R$-analogs of the symmetric algebras of a vector space $V$ and of its dual $V^*$. Let us precise that by a braiding
$R:\vv\to \vv$ we mean a solution of the braid relation
$$
R_{12}R_{23}R_{12} = R_{23}R_{12}R_{23},
$$
where $R_{12} = (R\ot I)$, $R_{23} = (I\ot R)$ and $I$ stands for the identity operator.

A fermionic version of a QD of Fock type can be defined in a similar way by using $R$-analogs $\Lambda_R(V)$ and $\Lambda_R(V^*)$ of the skew-symmetric algebras
$\Lambda(V)$ and $\Lambda(V^*)$.

In the classical case, i.e. when $R$ is the usual flip $P$,  the algebra $\mathrm{Sym}(V)$ often plays the role of a Fock space. Consequently, its generators give rise to creation
operators, whereas generators of the algebra  $\mathrm{Sym}(V^*)$ give rise to  annihilation ones. The permutation relations between the generators of these two algebras allows one 
to define an action of the annihilation operators on the algebra $\mathrm{Sym}(V)$ as partial derivatives. These two algebras form a classical version of a Fock double. By combining 
the annihilation and creation operators in a proper way we can perform a bosonic realization (or simply, a bosonization) of the Lie algebra $gl(N)$ and consequently its enveloping 
algebra $U(gl(N))$.

In the present paper we construct $R$-analogs  of these well-known objects, provided  braidings $R$ come from  the quantum groups (QG) $U_q(\ggg)$, corresponding to the
classical
simple Lie algebras $\ggg$. More precisely, the role of  such a braiding $R$ is played by the composition of the image of the universal $R$-matrix in $\vv$ and the usual flip
(here $V$
is the first fundamental space, called basic).

If $\ggg\in A_n$, then $R$ is a Hecke symmetry. Recall that by a Hecke symmetry we mean a braiding $R$ satisfying an additional relation
$$
(R-q\, I)(R+\qq I)=0,\quad q\not=\pm 1.
$$
If $\ggg$ belongs to other series of classical Lie algebras, then $R$ is a Birman-Murakami-Wenzl symmetry (see Section 4). For all corresponding QG $U_q(\ggg)$ the
$R$-symmetric
algebra $\mathrm{Sym}_R(V)$ and $R$-skew-symmetric algebras $\Lambda_R(V)$ of the space $V$ are well defined. Similar algebras for the dual space $V^*$ are also well defined.
Hopefully, all these algebras are deformations of their classical counterparts (see Remark \ref{FRT}). The first aim of the present paper is to construct the QD from the
above algebras.

The second aim consists in studying the following problem: which algebras can be constructed by means of  the bosonization (or the fermionization) procedure in a way similar
to the
classical one. It should be emphasized that the final algebras are not  the corresponding Quantum Groups themselves but they are covariant with respect to their
actions\footnote{In this connection we want to remark that there are numerous papers (we only mention the pioneering paper  \cite{M}), where the QG are realized via the
$q$-counterparts of harmonic oscillators.}. Moreover, the final algebras differ drastically from each other for different classical series. If $R$ comes from the QG
$U_q(sl(N))$ (in what follows we call it the standard symmetry), we
get the so-called modified reflection equation algebra $\LL(R)$, which is  generated by the unit\footnote{All algebras below are assumed to be unital.} and elements $l_i^j$,
subject to the following system of relations:
\be
R_{12}L_1 R_{12}L_1-L_1R_{12} L_1 R_{12} = R_{12}L_1 -  L_1 R_{12}, \qquad L=\|l_i^j\|_{1\leq i,j\leq N},\quad L_1=L\ot I.
\label{mRE}
\ee
This algebra is a deformation of the universal enveloping algebra $U(gl(N))$. Observe that the algebra $\LL(R)$ itself can be treated as an enveloping algebra of a
generalized
(or braided) Lie algebra defined in the space $\End(V)$ of internal endomorphisms. Moreover, a similar algebra $\LL(R)$ can be associated with any Hecke symmertry $R$ (or
with an involutive symmetry, i.e. such that $R^2=I$), provided $R$ is skew-invertible (see the next section for definition). A big family of such involutive and Hecke symmetries
has
been constructed in \cite{G}. For each of these symmetries the corresponding algebra $\LL(R)$ can be constructed by means of the bosonization or the fermionization procedures
and can be given a meaning of the enveloping algebra of a generalized Lie algebra. Moreover, $\LL(R)$ can be endowed with the trace $Tr_R$, which is coordinated with the
mentioned generalized Lie algebra structure.

Especially we are interested in the following question: to what extent the above construction, associated with Hecke symmetries $R$, can be generalised to the case, when $R$
is a BMW symmetry coming from the  orthogonal or simplectic QG. In this case we construct some QD associated with such symmetries and introduce algebras arising from the
bosonization or
fermionization procedure. It turns out that the final algebras differ from $\LL(R)$ and do not have any meaning of an enveloping algebras. This result shows that the BMW
symmetries are not well adapted to constructing  Lie-algebraic objects. By contrast, they fit well for introducing  group-like objects.
We refer the reader to the end of Section 4, where we discuss this phenomenon.

Besides, we apply a similar scheme to the case when $R$ is a current  braiding constructed by means of the so-called Baxterization procedure from an involutive or Hecke
symmetry.
We construct the QD of Fock type in the spirit of the Zamolodchikov-Faddeev algebras and exhibit the algebra arising from the corresponding bosonic realization. It turns out
that this algebra is inhomogeneous counterpart of the braided Yangians, introduced in \cite{GS1}.

The paper is organized as follows. In  section 2 we recall some basic notions and constructions, related to braidings. In section 3 we introduce the corresponding QD of Fock
type and
perform the bosonic and fermionic realizations of the algebras $\LL(R)$. In section 4 we consider the algebras, associated with BMW symmetries coming from the QG of the
$B_n$,
$C_n$
and $D_n$ series. In section 5 we exhibit QD of Zamolodchikov--Faddeev type as well as the corresponding algebras obtained by the bosonization procedure.

\medskip

\noindent
{\bf Acknowledgements}
The work of P.S. was partially supported by RFBR grant no. 19-01-00726.

\section{Hecke symmetries and corresponding algebras}

Let $R:\vv\to \vv$ be a braiding. By a successive application of  $R$ it is possible to transpose arbitrary tensor powers $V^{\ot k}$ and $V^{\ot l}$, that is to define a
linear map
$$
V^{\ot k}\ot V^{\ot l} \to V^{\ot l}\ot V^{\ot k}, \qquad k,l\in \Bbb{N}.
$$
However, we also need a transposition rule of the space $V$ and its dual $V^*$. For this purpose we should extend the action of the operator $R$ on spaces $V\ot V^*$
and $V^*\ot V$. A method of defining such an extension belongs to V.Lyubashenko \cite{L1,L2}.

We fix a basis $\{x_i\}_{1\leq i \leq N}$ of the space $V$ and  represent  the operator $R$ in the corresponding basis of the space $V^{\ot 2}$ by an $N^2\times N^2$ matrix
$\|R_{i j}^{k l}\|$:
$$
R(x_i\ot x_j)=R_{i j}^{k l}\,x_k\ot x_l.
$$
A summation over the repeated indices is always uderstood.

Let us fix a nondegenerate bilinear form $<\,,\,>_r: \,V\otimes V^*\rightarrow {\Bbb C}$ and choose {\it a right dual} basis\footnote{The term ``right'' reflects the position
of the argument $x^j$ in the bilinear form $<\,,\,>_r$.} $\{x^i\}_{1\le i\le N}$ of $V^*$ with respect to this form:
$<x_i,x^j>_r = \delta_i^j$. The pairing of the spaces $V^{\ot 2}$ and  $(V^*)^{\ot 2}$ is defined by
$$
<x_i\ot x_j, x^k\ot x^l>_r=<x_i,x^l>_r<x_j,x^k>_r.
$$

The operator $R$ on the space $(V^*)^{\ot 2}$ is defined as an operator adjoint to $R$, and its action on the right dual basis vectors is given by the rule:
$$
R(x^i\ot x^j)=R^{ji}_{lk}\,x^k\ot x^l.
$$
We keep the notation $R$ for this adjoint operator.

A braiding  $R$ is called {\em skew-invertible} if  there exists an operator $\Psi:\vv\to \vv$ such that
$$
\mathrm{ Tr}_{2} R_{12}\Psi_{23}=P_{13} \quad \Leftrightarrow\quad R_{ij}^{kl}\, \Psi_{lm}^{jn}=\de_m^k\, \de_i^n.
$$
As usual, the bottom indexes in the notation like $\Psi_{23}$ indicate the position of components of a tensor product where the operator acts.
The aforementioned extension of the braiding $R$ on the spaces $V\ot V^*$ and $V^*\ot V$ reads (see \cite{GPS1} for detail):
$$
R(x_i\otimes x^j) = x^k\ot x_l (R^{-1})^{lj}_{ki}, \qquad R(x^i\otimes x_j) = x_l\ot x^k \,\Psi^{li}_{kj}.
$$
If $R$ is skew-invertible we can define the following operators (matrices)
\be
B=\mathrm{Tr}_{1} \Psi_{12}\,\, \Leftrightarrow \,\,B_i^j= \Psi_{ki}^{kj}, \qquad C=\mathrm{Tr}_{2} \Psi_{12}  \,\,\Leftrightarrow \,\,C_i^j= \Psi_{ik}^{jk}.
\label{BC}
\ee
If the matrices $B$ and $C$ are nonsingular, the symmetry $R$ is called {\it strictly} skew-invertible.

Introduce now the second bilinear form $<\,,\,>_l: V^*\ot V\rightarrow {\Bbb C}$ as the composition of the linear maps:
$$
<\,,\,>_l  \,= \, <\,,\,>_r\circ R.
$$
Then we can introduce {\it a left dual} basis $\{\tilde x^i\}_{1\le i\le N}$ of $V^*$ by the requirement
$$
<\tilde x^i,x_j>_l = \delta_j^i,\quad 1\le \forall\,i,j\le N.
$$
If $R$ is strictly skew-invertible, then it is not difficult to establish the following relations
\be
<x^j, x_i>_l\,=\,B_i^j,\qquad <x_i, \tilde x^j>_r\,=\,(B^{-1})_i^j.
\label{BBCC}
\ee
Note that in the standard case the pairings $<\,,\,>_r$ and $<\,,\,>_l$ are $U_q(sl(N))$-covariant.

Now, consider the space of the internal endomorphisms
\be
\End(V)\cong V\ot V^*
\label{ident}
\ee
equipped with the basis $l_i^j=x_i\, x^j$. Hereafter, we  omit the sign $\ot $, if it does not lead to misunderstanding.

By using the identification (\ref{ident}), we get the following multiplication table in this
basis
$$
l_i^j\circ l_k^m=B_k^j \, l_i^m.
$$
 By contrast with the matrix $B$, determining the product $\circ$, the matrix $C$ enters the definition of  the so-called  $R$-trace. Namely, for any $N\times N$ matrix $M$
 (even with noncommutative entries)  its $R$-trace is defined by the following formula $\mathrm{Tr}_R M=\mathrm{Tr}\, CM$.

Now, consider the following map
\be
\RR: \End(V)^{\ot 2}\to \End(V)^{\ot 2},\qquad
R_{12} L_1 R_{12} L_1\stackrel{\RR}{\mapsto} L_1 R_{12} L_1R_{12},
\label{R}
\ee
where $L=\|l_i^j\|_{1\leq i, j\leq N}$ and $L_1=L\ot I$. A more explicit expression of this map can be found in \cite{GPS1}.

Using the map $\RR$, we can define an analog of the Lie bracket in the space $\End(V)$ by putting
\be
[X,Y]_R=X\circ Y-\circ \RR(X\ot Y), \qquad X,Y\in \End(V). \label{brbr}
\ee
The  data $(\End(V), \RR, [\,\,,\,]_R)$ is called a {\em generalized (or braided) Lie algebra} and is denoted $gl(V_R)$.

Observe that for this $gl(V_R)$-bracket the following analog of the Jacobi identity is valid on $\End(V)^{\ot 3}$ (here we omit the subscript  $R$)
$$
[\,,\,]\,[\,,\,]_{23}(I-\RR)_{23}=[\,,\,]\,[\,,\,]_{12}.
$$
A proof of this relation is given in \cite{GPS1}.

Note that if $R$ is an involutive symmetry, the corresponding Jacobi identity can be cast in the form similar to the classical one:
$$
 [\,,\,][\,,\,]_{23}(I+\RR_{12}\,\RR_{23}+\RR_{23}\,\RR_{12})=0.
$$

Let us explain  the reason why  the algebra $\LL(R)$ plays the role of the enveloping algebra of the braided Lie algebra $gl(V_R)$.
Indeed,   the defining system of the algebra $\LL(R)$ can be cast under the following form
\be
l_i^j \, l_k^l-\RR(l_i^j \, l_k^l)=[l_i^j,\, l_k^l]_R,
\label{mmRE}
\ee
where $\RR$ is the operator introduced in (\ref{R}). Thus, by applying the product $\circ$ to the left hand side of the relation (\ref{mmRE}), where we treat the
generators  $l_i^j$ as  elements of $\End(V)$, we get the element $[l_i^j,\, l_k^l]_R$.

Now, we want to show that the $R$-trace $\Tr_R$ is coordinated with the bracket $[\,,\,]_R$ in the following sense
$$
\mathrm{Tr}_R\, [X,Y]_R=0, \quad \forall\, X,Y\in \End(V),
$$
which is similar to the classical relation for the trace of the commutator of two operators.
To this end we first compute ${\Tr}_R\, l_i^j$. Since $l_i^j\triangleright x_k=B_k^j\, x_i$, then applying the $R$-trace to the image of the element $l_i^j$ we get
${\Tr}_R\, l_i^j= {\Tr}\,B_i^k\, C_k^j=\al \,\de_i^j$, where $\al$ is a nontrivial numerical factor, whose value is determined by the bi-rank of $R$ (see \cite{GPS1} for detail).
Thus, taking $X=l_i^j$ and $Y=l_k^m$, we have to show ${\Tr}_R\, [l_i^j,l_k^m]_R=0$. This equality is a consequence of the following chain of relations:
$$
{\Tr}_R \circ(R_{12} L_1R_{12} L_1-L_1 R_{12} L_1 R_{12})={\Tr}_R(R_{12} L_1-L_1 R_{12})=\al(R_{12} I-I\, R_{12})=0.
$$

With any Hecke symmetry $R$ we associate the $R$-symmetric $\mathrm{Sym}_R(V)$ and $R$-skew-symmetric $\Lambda_R(V)$ algebras of the spaces
$V$ and $V^*$ by setting
\be
 \mathrm{Sym}_R(V)=T(V)/\langle\mathrm{Im}(q\, I-R)\rangle =T(V)/\langle q\, x_i\, x_j-R_{i\, j}^{k\, l}\,x_k\, x_l\rangle ,
\label{symV}
\ee
$$
\Lambda_R(V)=T(V)/\langle \mathrm{Im}(\qq\, I+R)\rangle =T(V)/\langle \qq x_i\, x_j+R_{i\, j}^{k\, l}\,x_k\, x_l\rangle ,
$$
\be
\mathrm{Sym}_R(V^*)=T(V^*)/\langle \mathrm{Im}(q\, I-R)\rangle =T(V^*)/\langle q\, x^j\, x^i-R^{i\, j}_{k\, l}\,x^l\, x^k\rangle ,
\label{symV*}
\ee
$$
\Lambda_R(V^*)=T(V^*)/\langle \mathrm{Im}(\qq\, I+R)\rangle =T(V^*)/\langle \qq\, x^j\, x^i+R^{i\, j}_{k\, l}\,x^l\, x^k\rangle,
$$
where
$$
T(V)=\bigoplus_k V^{\ot k},\quad T(V^*)=\bigoplus_k V^{*\,\ot k}
$$
are the free tensor algebras and $\langle J\,\rangle$ stands for a two-sided ideal generated by a subset $J$ in the given algebra. If $R$ is an involutive symmetry,
we put $q=1$ in all these formulae.

It should be emphasized that all these algebras are deformations of their classical counterparts, provided $R$ is a deformation of the usual flip $P$, i.e. dimensions of
their
homogenous components are classical (for a generic $q$ if $R$ is Hecke). A similar statement is valid, if a Hecke symmetry is a deformation of another involutive symmetry,
for instance, a super-symmetry.

To complete this section, we also mention two quantum matrix algebras, associated with any braiding $R$. One of them is the so-called RTT algebra, defined by the system
\be
R_{12} T_1\, T_2=T_1\, T_2R_{12},\qquad T=\|t_i^j\|_{1\leq i, j \leq N}, \quad T_1=T\ot I, \quad T_2=I\ot T.\label{RTT}
\ee
The other one is the RE algebra defined via (\ref{mRE}) but with the vanishing right hand side. We do not use the former algebra. The latter one appears only under its
modified form. However, it is worth noting that  if $R$ is a Hecke symmetry,  the RE algebra and its modified version are isomorphic to each other.

\section{Quantum doubles and bosonic (fermionic) realization of $\LL(R)$}

Our next aim is to construct some quantum doubles (QD) from the above $R$-symmetric and $R$-skew-symmetric algebras and use them in order to perform  bosonic
and fermionic realizations of the algebra $\LL(R)$.

First, we recall the definition of doubles of associative algebras from \cite{GS3}. Consider two associative unital algebras $A$ and $B$ equipped with a linear map
$$
\sigma :A \ot B\to B\ot A,
$$
which satisfies the following conditions:
$$
\sigma\circ (\mu_A\ot \id)=(\id\ot \mu_A)\circ \sigma_{12}\sigma_{23}\quad {\rm on}\,\, A\ot  A\ot B,
$$
$$
\sigma\circ (\id\ot \mu_B)=(\mu_B\ot \id)\circ \sigma_{23}\sigma_{12}\quad {\rm on}\,\, A\ot  B\ot B,
$$
$$
\sigma(1_A\ot b)=b\ot 1_A,\quad \sigma (a\ot 1_B)=1_B\ot a,\quad \forall\, a\in A, \,\forall\,b\in B,
$$
where $\mu_A:A\ot A\to A$ is the multiplication operation in the algebra $A$, $1_A$ is the unit element of $A$ and the symbols $\mu_B$ and $1_B$ have the same meaning
in the algebra $B$. The symbol ``$\id$'' stands for the identity operator.

It is not difficult to see that the bilinear map $*:(B\ot A)^{\ot 2} \rightarrow B\ot A$ defined by the rule
$$
(b\ot a) * (b'\ot a'):=(\mu_B\otimes \mu_A)\circ (\mathrm{id}_B\otimes\sigma_{23}\otimes\mathrm{id}_A) (b\ot a\ot b'\ot a')
$$
endows  the space $B\ot A$ with the structure of a unital associative algebra and its unit element is $1_B\ot 1_A$.

\begin{definition}\rm
We call the data  $(A,\,B,\,\sigma)$ and the corresponding  algebra $(B\ot A, *)$ a {\it double of associative algebras} and denote it $(A,B)$.

 An associative double $(A,B)$ is called {\it quantum} if the map $\sigma$ is constructed with the use of a braiding, different from a (super-)flip.
 \end{definition}

Observe that if the both algebras $A$ and $B$ are introduced via some systems of relations on the generators, the above conditions on the data $(A,\,B,\,\sigma)$ mean
that the permutation map  preserves these relations. In this sense we say that the permutation map and the defining relations of the algebras $A$ and $B$ are {\em compatible}.

Now, assume  the algebra $A$ to be equipped with a counit (an algebra homomorphism)\footnote{In general we do not consider any coalgebraic structure in the algebra $A$.
So, the counit is only coordinated with the algebraic structure.} $\varepsilon_A:A\to \C$. Then it is possible to define a linear action  of the algebra $A$ on $B$ by setting
$$
a\triangleright b = (\mathrm{id}_B\otimes\varepsilon_A)\circ \sigma(a\ot b),\quad \forall\,a\in A,\, \forall\,b \in B.
$$
Identifying $b\ot 1_\C$ and $b$, we get that $a\triangleright b\in B$. Thus,  each element $a\in A$ defines a linear operator
$$
Op(a):B\to B.
$$
It is easy to check that the map $a\mapsto Op(a)$ defines a representation of the algebra $A$ in the algebra $B$:
$$
Op(ab)=Op(a) Op(b),\quad Op(1_A)=\mathrm{id}_B.
$$

Besides, the aforementioned  Heisenberg double, we can also mention QD  composed from two copies of the RTT algebras. The permutation relations between them are also
similar to the system (\ref{RTT}). A counit $\varepsilon$ defined on one component of such a double is usually defined by $\varepsilon(t_i^j)=\de_i^j$.
Note that the corresponding representation of the   component $A$, defined according to the above scheme, is  a deformation of the trivial one
$\pi(T)=I$ and consequently is not interesting indeed. The QD of Fock type, which we are dealing with, are more meaningful. We introduce these doubles, assuming $R$ to be
a Hecke symmetry.

Let $R:\vv\to\vv$ be a skew-invertible Hecke symmetry, $A=\mathrm{Sym}_R(V^*)$ and $B=\mathrm{Sym}_R(V)$. Construct a quantum double $(A,B)$ by
introducing the following permutation relations\footnote{By permutation relations we mean the equality of the form $a\ot b=\sigma(a\ot b)$.}:
\be
x^a R_{a\, i}^{b\,j}\, x_b=\qq \, x_i\, x^j+\de_i^j \quad \Leftrightarrow \quad
x^l\, x_k =\qq\, x_i\, x^j\, \Psi_{j \,k}^{\,i\, l}+B_k^l.
\label{perm}
\ee
\begin{remark} \rm Recall, that we use the right dual basis in the space $V^*$. For the left dual basis $\{\tilde{x}^j\}_{1\le j\le N}$, the permutation relations read
$$
x_b\, R_{i\,a}^{j\,b}\, \tx^a=\qq \,\tx^j\,  x_i+\de_i^j \quad\Leftrightarrow\quad
x_k\, \tx^l =\qq\, \tx^j\, x_i\,  \Psi_{k \,j}^{\,l\,i}+C_k^l.
$$
\end{remark}

\begin{proposition} \label{prop:Hecke}
The defining systems of the algebras $\mathrm{Sym}_R(V)$ and $\mathrm{Sym}_R(V^*)$ are compatible with the permutation relations
{\rm (\ref{perm})}.
\end{proposition}

\noindent
{\bf Proof.} Below, we use the Dirac's ``bra''  and ``ket''  notation:
$$
x_i\rightarrow x_{|1\rangle},\quad x_ix_j\rightarrow x_{|1\rangle}x_{|2\rangle},\quad  x^i\rightarrow x^{\langle 1|},\quad x^ix^j\rightarrow x^{\langle 1|}x^{\langle
2|},\,\,\dots
$$
We call this form of  notation the {\it matrix} one. Some typical examples are
$$
x^{\langle 1|}x_{|1\rangle} = \sum_{i=1}^N x^ix_i,\qquad x_{|1\rangle}x^{\langle 1|} = x_ix^j,   \qquad
 x^{\langle 2|}R_{|12\rangle}^{\,\langle 12|}x_{|2\rangle} =
\sum_{a,b=1}^Nx^a R_{\,i\,a}^{\,\,j\,b}x_b.
$$
In what follows the notation $R_{|12\rangle}^{\,\langle 12|}$ will be simplified  to $R_{12}$.

Thus,  the defining relations of the algebras (\ref{symV}), (\ref{symV*}) and the permutation relations (\ref{perm}) can be respectively rewritten as:
$$
\begin{array}{lcl}
R_{\,i\, j}^{\,\,k\, l}\,x_k\, x_l = q\, x_i\, x_j &\quad \rightarrow\quad&R_{12}x_{|1\rangle}x_{|2\rangle} = q \,x_{|1\rangle}x_{|2\rangle}\\
\rule{0pt}{6mm}
x^l\, x^kR^{\,\,i\, j}_{\,k\, l} = q\, x^j\, x^i &\rightarrow&x^{\langle 2|}x^{\langle1|}R_{12} = q\,x^{\langle 2|}x^{\langle1|}\\
\rule{0pt}{6mm}
x^a R_{\,a \,i}^{\,\,b\,j}\, x_b=\qq \, x_i\, x^j+\de_i^j  &\rightarrow & x^{\langle 1|}R_{12}x_{|1\rangle} =
\qq\,x_{|2\rangle}x^{\langle 2|} + I_{|2\rangle}^{\,\langle 2|}\\
\rule{0pt}{6mm}
&\mathrm{or}& \,x_{|2\rangle}x^{\langle 2|} = q\,x^{\langle 1|}R_{12}x_{|1\rangle} - q\,I_{|2\rangle}^{\,\langle 2|}.
\end{array}
$$

Let us verify, for example, that the defining relations of the algebra $\mathrm{Sym}_R(V)$ (\ref{symV}) are not destroyed by the permutation with vectors from $V^*$.
For this purpose we permute an arbitrary basis vector $x^i$ with elements of  the generating set of the ideal in (\ref{symV}):
\begin{eqnarray*}
( R_{23} \!\!\!\!&-&\!\!\!\! qI_{23})x_{|2\rangle}x_{|3\rangle} x^{\langle 3|} =( R_{23} - qI_{23})x_{|2\rangle}
\left( q\,x^{\langle 2|}R_{23} x_{|2\rangle} - q\,I_{|3\rangle}^{\,\langle 3|} \right)\\
&=&\!\!\! -q ( R_{23} - qI_{23})x_{|2\rangle}+q( R_{23} - qI_{23})\Big(qx^{\langle 1|}R_{12}x_{|1\rangle} -q I_{|2\rangle}^{\langle 2|}\Big)R_{23}x_{|2\rangle}\\
&=& q^2x^{\langle 1|}(R_{23} - qI_{23})R_{12}R_{23}x_{|1\rangle}x_{|2\rangle} = q^2x^{\langle 1|}R_{12}R_{23}(R_{12} - qI_{12})x_{|1\rangle}x_{|2\rangle} .
\end{eqnarray*}
Here, by passing to the third line we have taken into account the consequence of the Hecke condition $(R_{23} - qI_{23})R_{23} = -\qq (R_{23} - qI_{23})$,
while the last equality is due to the braid relation for $R$.

So, we finally get the equality:
$$
( R_{23} - qI_{23})x_{|2\rangle}x_{|3\rangle}\, x^{\langle 3|} = q^2x^{\langle 1|}R_{12}R_{23}\,(R_{12} - qI_{12})x_{|1\rangle}x_{|2\rangle},
$$
which means the compatibility of the defining relations of the algebra $\mathrm{Sym}_R(V)$ with the permutation relations (\ref{perm}).

The compatibility condition for the defining relations of  $\mathrm{Sym}_R(V^*)$ and the permutation relations can be verified in the same way.\hfill\rule{6.5pt}{6.5pt}

\medskip

Now, we  introduce the counit similarly to the classical case:
$$
\varepsilon_A:A=\mathrm{Sym}_R(V^*)\to \C,\qquad \varepsilon_A(1_A)=1_{\C},\quad \varepsilon_A(x^j)=0,\quad 1\le \forall j\le N.
$$

According to the above scheme all elements of the algebra $A=\mathrm{Sym}_R(V^*)$ can be represented by  operators acting on the algebra $B=\mathrm{Sym}_R(V)$.
Compute the action of elements $x^j$ onto some monomials in $x_i$:
\be
x^j \triangleright 1_B=0,\quad   x^j \triangleright x_i=B_i^j,\quad x^j \triangleright (x_i\, x_k)=B_i^j \, x_k+\qq\, \, B^l_k\Psi_{li}^{mj}\, x_m,\quad {\rm etc}.
\label{act}
\ee
Note that using the properties of $\Psi $ and $R$ we can rewrite the above actions in the {\it equivalent} matix form:
$$
x^{\langle 1|}\triangleright\Big(R_{12}x_{|1\rangle}\Big) = I_{|2\rangle}^{\,\langle 2|}, \qquad
x^{\langle 1|}\triangleright\Big(R_{12}R_{23}x_{|1\rangle}x_{|2\rangle}\Big) = (R_{23} +\qq I_{23})x_{|2\rangle}.
$$
As we noticed above, the generators of the  algebra $A=\mathrm{Sym}_R(V^*)$ are analogs of  the usual annihilation operators, whereas the generators of
$B=\mathrm{Sym}_R(V)$ are analogs of  the creation operators: the result of the action of $x_i$ on  $a\in \mathrm{Sym}_R(V)$ is the product $x_i a$. The above
double $(A, B)$ is said to be {\em quantum bosonic double of  Fock type}.

In a similar manner it is possible to construct a {\em quantum fermionic double of Fock type} by setting $A={\mathrm{\Lambda}}_R(V^*)$, $B={\mathrm{\Lambda}}_R(V)$
and
\be
x^a R_{a\, i}^{\,b\,j}\, x_b=-q \, x_i\, x^j+\de_i^j \quad \Leftrightarrow\quad
x^l\, x_k =-q\, x_i\, x^j\, \Psi_{j \,k}^{\,i\, l}+B_k^l.
\label{perm1}
\ee
In this case all defining relations are also compatible.

Now, compose the matrix $L=\|l_i^j\|$, $l_i^j=x_i x^j$, where $x_i\in B$, $x^j\in A$ and $(A, B)$ is a bosonic or fermionic double of Fock type.
Then the following proposition is valid.

\begin{proposition}
The matrix $L$ meets the relations {\rm (\ref{mRE})}.
\end{proposition}

\noindent
{\bf Proof.} In the bosonic case this proposition was proved in \cite{GS1}. In the fermionic case the proof is similar.\hfill \rule{6.5pt}{6.5pt}

\medskip

Thus, representing the factors of  the product $x_i x^j$  as the creation and annihilation operators acting on the algebra $B$, we get the so-called {\em bosonic (or
fermionic) realization} of the algebra $\LL(R)$ similar to that of the algebra\footnote{Note that if $R\to P$ as $q\to 1$ (for instance, if $R$ is a standard symmetry) the system
(\ref{mRE}) turns into the matrix equality
$$
P L_1 P L_1-L_1 P L_1 P=P L_1- L_1 P\quad\Leftrightarrow  \quad   l_i^j\,l_k^m-l_k^m\,l_i^j= l_i^m\,\de_k^j-l_k^j\,\de_i^m,
$$
valid for the usual generators of the algebra $U(gl(N))$.} $U(gl(N))$.

Moreover, similarly to the classical case  the generators $l_i^j$ preserve any homogenous component $\mathrm{Sym}_R^{k}(V)$, $k\ge 1$. Thus, we get a
series of finite-dimensional representations of the algebra $\LL(R)$. For instance,  the action of the elements $l_i^j=x_ix^j$ onto the first component
$\mathrm{Sym}_R^{1}(V)=V$ is as follows
\be
l_i^j\triangleright x_k=x_i\triangleright (x^j\triangleright x_k)= B^j_k\, x_i.
\label{acttt}
\ee

\begin{remark}\rm
\label{remm}
Note that on the role of the algebra $B$ we can assign the free tensor algebra $T(V)$ instead of $\mathrm{Sym}_R(V)$ and keep the same permutation relations.
By considering all invariant subspaces with respect to the action $\triangleright$, it is possible to construct a more rich representation category of the algebra $\LL(R)$.
Such a category was constructed \cite{GPS1} by means
of  a ``braided bi-algebra structure" of $\LL(R)$. Our present method based on bosonization  is a sense close to the method of constructing
differential calculus on the quantum hyperplane from \cite{WZ, H}.
\end{remark}

\section{Quantum doubles associated with BMW symmetries}

Let us recall that by a BMW symmetry we mean a braiding $R$ which comes from a QG belonging to one of the orthogonal ($B_n,\, D_n$) or simplectic ($C_n$) series. An explicit
form of these symmetries are exhibited in \cite{FRT}. Each of them is subject to some relations. We need only one of them:
\be
 (R-q  I)(R+\qq  I)(R-\mu I )=0.
\label{degree}
\ee
Note that if $R$ comes from the QG corresponding to an orthogonal group  then  $\mu=q^{1-N}$ and if $R$ comes from the QG corresponding to a simplectic group, then
$\mu=-q^{-1-N}$. Here $N$ is the dimension of the basic vector space.  Below, for the sake of simplicity we speak about orthogonal and simplectic QG.

Also, there are known examples of the BMW symmetries coming from super-QG with other values of $\mu$ (see \cite{I}).   Below, we only deal with the BMW symmetries, coming
from orthogonal or simplectic  QG. Consequently, the
parameter $\mu$ is assumed to take one of the above values. Thus, as $q\to 1$, this parameter tends to 1 for the orthogonal QG and to -1 for the simplectic ones.

In \cite{OP} there were considered the so-called BMW algebras, which are defined by similar relations but with the parameter $\mu$ independent on $q$.

In virtue of (\ref{degree}) there exist three complementary idempotents ${\PP}^q$, ${\PP}^{-\qq}$ and ${\PP}^{\mu}$ such that
$$
R=q\, {\PP}^q-\qq {\PP}^{-\qq}+\mu {\PP}^{\mu}.
$$
Each of them is a projector onto the corresponding eigenspace and can be explicitly written in terms of $R$. For instance, the projector
\be
{\PP}^{-q^{-1}}=\frac{(R-q I)(R-\mu I)}{(q+\qq)(\qq+\mu)}
\label{projP}
\ee
maps the space $\vv$ onto the eigenspace of the operator $R$ corresponding to the eigenvalue $-\qq$. Note that for a generic $q$ the denominator of this formula
does not vanish.

Now, let us define the $R$-symmetric and $R$-skew-symmetric algebras of the space $V$ by setting
\be
\mathrm{Sym}_R(V)=T(V)/\langle\mathrm{Im}\, {\PP}^{-q^{-1}}\rangle,\qquad {\mathrm{\Lambda}}_R(V)=T(V)/\langle\mathrm{Ker}\, {\PP}^{-q^{-1}} \rangle,
\label{sV}
\ee
provided $R$ comes from an orthogonal QG  and
$$
\mathrm{Sym}_R(V)=T(V)/\langle\mathrm{Ker} \, {\PP}^{q}\rangle,\qquad {\mathrm{\Lambda}}_R(V)=T(V)/\langle\mathrm{Im}\, {\PP}^{q}\rangle,
$$
provided $R$ comes from a simplectic QG.

\begin{remark} \label{FRT} \rm 
Note that according to \cite{FRT} $R$-symmetric algebras have the classical Poincar\'{e} series, i.e. these series are equal to the series
corresponding to the algebras $\mathrm{Sym}(V)$. However, any prof of this claim has never been published. We do not use this property. We want only
to observe that the space $\mathrm{Ker} \,{\PP}^{\mu}$ belongs to the second homogenous $R$-symmetric component of the space $\vv$ if the 
QG is  orthogonal  and to $R$-skew-symmetric component if QG is  simplectic.
\end{remark}

In a similar manner we define the $R$-symmetric and $R$-skew-symmetric algebras of the space $V^*$. For instance, we put (in the orthogonal case)
\be
\mathrm{Sym}_R(V^*)=T(V^*)/\langle\mathrm{Im}\, {\PP}^{-q^{-1}}\rangle,
\label{sV*}
\ee
where the idempotents ${\PP}^q$, ${\PP}^{-\qq}$ and ${\PP}^{\mu}$ acting in the space $(V^*)^{\ot 2}$ are expressed in terms of $R$ by the same formulae as above. 
We only have to take into account the way of extending the braiding $R$ onto the space $V^{*\otimes 2}$, namely,  ${R}(x^k\ot x^l)={R}_{ji}^{lk}(x^i\ot x^j)$.

Let us define the permutation relations by formula (\ref{perm}) for an orthogonal QG and by (\ref{perm1}) for a simplectic QG. Then the following claim is valid.

\begin{proposition}
 \label{odin}
The permutation relations are compatible with the defining systems of the algebras $\mathrm{Sym}_R(V)$ and $\mathrm{Sym}_R(V^*)$ for the orthogonal series
and with these of the algebras ${\mathrm{\Lambda}}_R(V)$ and ${\mathrm{\Lambda}}_R(V^*)$ for the simplectic series.
\end{proposition}

\noindent
{\bf Proof.} The proposition is proved by a straightforward calculation similar to the case of the Hecke symmetries (see Proposition \ref{prop:Hecke} above).
Nevertheless, we give a short sketch of the proof in order to stress the difference with the Hecke case. 

Consider the algebra $\mathrm{Sym}_{R}(V)$ for the orthogonal QG defined by the quotient (\ref{sV}). Taking into account the permutation relations (\ref{perm}) we get the following transformation:
$$
x_{|2\rangle}x_{|3\rangle}x^{\langle 3|} = q^2x^{\langle 1|}R_{12}R_{23}x_{|1\rangle}x_{|2\rangle} -q^2(R_{23} +q^{-1}I_{23})\,x_{|2\rangle}.
$$  
Then for the generating set of the ideal $\langle\,\mathrm{Im}\, {\PP}^{-q^{-1}}\rangle$ we have:
$$
{\PP}^{-q^{-1}}(R_{23})x_{|2\rangle}x_{|3\rangle}x^{\langle 3|}  = q^2x^{\langle 1|}R_{12}R_{23}{\PP}^{-q^{-1}}(R_{12})x_{|1\rangle}x_{|2\rangle}  -
q^2{\PP}^{-q^{-1}}(R_{23})(R_{23} +q^{-1}I_{23})\,x_{|2\rangle}.
$$ 
The first therm in the right hand side is a consequence of the braid relation on $R$:
$$
{\PP}^{-q^{-1}}(R_{23})R_{12}R_{23} = R_{12}R_{23}{\PP}^{-q^{-1}}(R_{12}),
$$
while the second term is equal to zero due to the definition (\ref{projP}) of the projector ${\PP}^{-q^{-1}}$ and the cubic minimal polynomial (\ref{degree}) 
of the symmetry $R$. So, we finally get
$$
\underline{{\PP}^{-q^{-1}}(R_{23})x_{|2\rangle}x_{|3\rangle}}x^{\langle 3|}  = q^2x^{\langle 1|}R_{12}R_{23}\,\underline{{\PP}^{-q^{-1}}(R_{12})x_{|1\rangle}x_{|2\rangle}},
$$
which means the compatibility of the permutatuon relations (\ref{perm}) with the defining relations of the algebra $\mathrm{Sym}_R(V)$ for the
orthogonal case. All other compatibilities are verified by the similar calculations.\hfill \rule{6.5pt}{6.5pt}

\medskip

Now, consider the  QD $(A,B)$ of Fock type, corresponding to a BMW symmetry $R$. For $R$ corresponding to the orthogonal QG
we set $A=\mathrm{Sym}_R(V^*)$, $B=\mathrm{Sym}_R(V)$, and for $R$ corresponding to a simplectic QG we set  $A={\Lambda}_R(V^*)$, $B={\Lambda}_R(V)$.
By introducing the same counit as in the algebra $\LL(R)$, we get an action of the algebra $A$ on the algebra $B$.  The creation and annihilation operators are similar to
these from the bosonic (resp., fermionic) realization of the algebra $\LL(R)$.

Now, we introduce a matrix $L=\|l_i^j\|$ with entries $l_i^j=x_i x^j$.

\begin{proposition}
\label{dva}
For all series the matrix $L$ meets the following system
{\rm
\be
{\PP}_{12}L_1 R_{12} L_1-L_1R_{12} L_1{\PP}_{12}={\PP}_{12} L_1-L_1{\PP}_{12},
\label{new}
\ee
}
where ${\PP}={\PP}^{q}+{\PP}^{\mu}$ for the orthogonal series and ${\PP}={\PP}^{-\qq}+{\PP}^{\mu}$ for the simplectic one.
\end{proposition}

\noindent
{\bf Proof.}  Consider the case of $R$ corresponding to an orthogonal QG in detail. For a simplectic QG all considerations are analogous.
Below, we use the matrix notation introduced above.

Since the idempotents ${\PP}^q$, ${\PP}^{-\qq}$ and ${\PP}^{\mu}$ are complementary, we have
$$
 {\cal P}^{-\qq} + {\cal P}^{q} + {\cal P}^{\mu} = I.
 $$
Also, according to the definitions of the algebras $\mathrm{Sym}_R(V)$ and $\mathrm{Sym}_R(V^*)$ and the permutation relations the complete system on  the generators $x_i$
and $x^i$ are
$$
{\cal P}_{12}x_{|1\rangle}x_{|2\rangle} = x_{|1\rangle}x_{|2\rangle}, \qquad x^{\langle 2|}x^{\langle 1|}{\cal P}_{12} = x^{\langle 2|}x^{\langle 1|},
\qquad x^{\langle 1|}R_{12}x_{|1\rangle } = q^{-1}x_{|2\rangle}x^{\langle 2|}+I_{2},
$$
where ${\cal P} = {\cal P}^{q}+{\cal P}^{\mu}$.

Now, substitute  the matrix $L_1 = x_{|1\rangle}x^{\langle 1|}$ into the left hand side of (\ref{new}). Using the relations on the generators we transform the first summand
to the following expression:
$$
{\cal P}_{12}x_{|1\rangle}\underline{x^{\langle 1|}R_{12}x_{|1\rangle}}x^{\langle1|} = q^{-1}\underline{{\cal P}_{12}x_{|1\rangle}x_{|2\rangle}}x^{\langle 2|}
x^{\langle 1|} +{\cal P}_{12}x_{|1\rangle}x^{\langle1|} = q^{-1}x_{|1\rangle}x_{|2\rangle}x^{\langle 2|} x^{\langle 1|} +{\cal P}_{12}L_1.
$$
Here, underlined are the terms which undergo changes.
By the same steps of transformations the second summand can be rewritten in the form:
$$
x_{|1\rangle}\underline{x^{\langle 1|}R_{12}x_{|1\rangle}}x^{\langle1|}{\cal P}_{12} = q^{-1} x_{|1\rangle}x_{|2\rangle}\underline{x^{\langle 2|} x^{\langle 1|}{\cal P}_{12}}
+ x_{|1\rangle}x^{\langle1|}{\cal P}_{12} = q^{-1}x_{|1\rangle}x_{|2\rangle}x^{\langle 2|} x^{\langle 1|} + L_1{\cal P}_{12}.
$$
By taking the difference of the above expressions, we get (\ref{new}).\hfill \rule{6.5pt}{6.5pt}

\medskip

Observe that the algebra defined by (\ref{new})  cannot be treated as  the enveloping algebra $U(\ggg)$ of a generalized  Lie algebra, since the relation (\ref{new}) does not
enable us to define any operator $\rm{End}(V)^{\ot 2}\to \rm{End}(V)^{\ot 2}$
similar to $\RR$\footnote{Observe that if $R$ is a Hecke symmetry, the corresponding modified RE algebra can be cast in the form (\ref{new}), but with ${\PP=\PP}^q$.}.

Moreover, though the algebras $\mathrm{Sym}_R(V)$ and $\Lambda_R(V)$  are well defined in this case  and are (hopefully!) deformations of their classical counterparts, 
the algebras obtained via the bosonic (fermionic) realization are not deformations of $U(\ggg)$.

Nevertheless, proceeding in a similar way as described above we can construct finite dimensional representations of the algebra defined by (\ref{new}) in homogenous 
components of the algebras $\mathrm{Sym}_R(V)$ or  ${\Lambda}_R(V)$ depending on the series.

Concluding this section we note that if $R$ is an involutive or Hecke symmetry, deforming the usual flip $P$, the both algebras $\mathrm{Fun}(GL(N))$ and $U(gl(N))$ admit quantum
deformations. Moreover, there are two different deformations of the algebra $\mathrm{Fun}(GL(N))$: the corresponding RTT algebra or the RE one. However, a deformation 
of the algebra $U(gl(N))$ can be performed  only via the modified RE algebra.

If $\ggg$ is an orthogonal or simplectic Lie algebra, the space $\mathrm{Fun}(G)$  of functions on the corresponding group $G$ can be deformed as was done in  \cite{OP}. 
Note that the authors of the cited paper use Quantum Matrix Algebras of a general form including the RTT and RE algebras. By contrary, the algebra $U(\ggg)$ does not 
have any  deformation,  covariant with respect to the adjoint action of the corresponding QG $U_q(\ggg)$.

As we have noticed in Introduction,  the BMW symmetries are not well adapted to constructing $R$-analogs of Lie algebras.

\begin{remark} \rm Now, we would like to discuss the Poisson counterparts of the algebras obtained by the bosonization. If $R$ is a Hecke symmetry deforming the usual 
flip $P$ (for instance, that coming from the QG $U_q(sl(N))$ or the Crammer-Gervais one), then the Poisson counterpart of the algebra $\LL(R)$ (more precisely, of the 
algebra, obtained by introducing a second deformation parameter  $\h$ in the front of the right hand side of (\ref{mRE})) is a pencil defined on the commutative algebra 
$\mathrm{Sym}(gl(N))$. This pencil is generated by the linear Poisson-Lie bracket, associated with the Lie algebra $gl(N)$, and by a quadratic bracket, corresponding to the case 
$\h=0$. It is interesting to observe that almost all brackets of this pencil are not unimodular. Nevertheless, the quantum counterpart of this pencil (namely, the algebra
$\LL(R)$) is endowed with a trace, but this trace is not usual.

If $R$ is a BMW symmetry coming from an orthogonal or simplectic QG, a similar Poisson pencil does not exist. However, a Poisson bracket, corresponding to the 
deformation of the algebra $\mathrm{Fun}(G)$ can be defined.
\end{remark}

\section{Quantum double of Zamolodchikov-Faddeev type}

In this section we deal with  algebras associated with current braidings.
By a current braiding we mean an operator $R(u,v)$ depending on  spectral parameters subject to the relation
$$
R_{12}(u,v)\,R_{23}(u,w)\, R_{12}(v,w)=R_{23}(v,w)\, R_{12}(u,w)\, R_{23}(u,v).
$$
The current braidings, we are dealing with, arise from the Baxterization procedure applied to an involutive (resp., Hecke) symmetry. These current braiding are of the form
\be
R(u,v)=R-\frac{I}{u-v},\quad \mathrm{ resp.,} \quad  R(u,v)=R-\frac{(q-\qq)\, u\, I}{u-v},   \label{braidd}
\ee
where $R$ is any skew-invertible involutive (resp., Hecke) symmetry.  The former braiding $R(u,v)$ is called {\em rational}, the  latter one is
called {\em trigonometric}. The reader is referred to \cite{GS1} for details.

Let us pass to the normalized braidings $\RRR(u, v)=g(u,v)^{-1} R(u,v)$, where
$$
g(u,v)= 1-\frac{1}{u-v},\quad \mathrm{ resp.}\quad g(u,v)=q-\frac{(q-\qq)u}{u-v}
$$
provided $R(u,v)$ is a rational (resp., trigonometrical) braiding (\ref{braidd}).

It should be emphasized that these operators are still subject to the braid relation (\ref{braidd}) and are involutive in the following sense
$$
\RRR(u,v)\, \RRR(v,u)=I.
$$
In virtue of this relation the operators $\RRR(u,v)$ have two eigenvalues $\pm 1$. So, we  introduce an analog of the algebra
$\mathrm{Sym}_R(V)$  by the following relations:
\be
\RRR_{ij}^{kl}(u,v)x_k(u)\, x_l(v)  = x_i(v)\, x_j(u) \quad\Leftrightarrow \quad R_{ij}^{kl}(u,v)x_k(u)\, x_l(v) = g(u, v)\,x_i(v)\, x_j(u).
\label{chet}
\ee

Write this system in a more detailed form  for a trigonometric braiding
$$
q\, x_i(v)\,x_j(u)-R_{ij}^{kl}x_k(u)\,x_l(v)=\frac{(q-\qq)\, u}{u-v}(x_i(v)\,x_j(u)-x_i(u)\,x_j(v)).
$$
(The defining relations of the algebra $\Lambda_R(V)$ can be written in a similar way.)

However, if we try to represent the currents $ x_i(u)$ as a formal series in $u$ with Fourier coefficients $x_i[m]$ as follows
$$
x_i(u)=\sum_{m\in \Z} x_i[m] \, u^{-m-1},
$$
we find that it is not possible to rewrite the defining system via polynomial relations on these coefficients.
By following the classical pattern,  we  present the algebra
in terms of the so-called {\it half-currents}
$$
x_i^+(u)=\sum_{m\in \Z,\, m<0} x_i[m] \, u^{-m-1},\quad x_i^-(u)=\sum_{m\in \Z,\, m\geq 0} x_i[m] \, u^{-m-1},\quad x_i(u)=x_i^+(u)+x_i^-(u).
$$

Then by imposing  the relations
\be
\begin{array}{l}
g(u, v)\,x_i^{\pm}(v)\, x_j^{\pm}(u)=R_{ij}^{kl}(u,v)x_k^{\pm}(u)\, x_l^{\pm}(v),\\
\rule{0pt}{7mm}
g(u, v)\,x_i^-(v)\, x_j^+(u)=R_{ij}^{kl}(u,v)x_k^+(u)\, x_l^-(v),
\end{array}
\label{rela}
\ee
and by using the expansion $\frac{1}{u-v}=\sum_{p\geq 0}\frac{v^p}{u^{p+1}}$, we can rewrite (\ref{chet}) in terms of the Fourier coefficients.
However, the system  (\ref{chet}) for the currents $x_i(u)$ is preserved.

According to the scheme exhibited in Section 3, we introduce the  current $R$-symmetric algebra of the dual space by the following system
\be
g(u, v)\,x^i(u)\, x^j(v)=R^{ji}_{lk}(u,v)x^k(v)\, x^l(u).
\label{rela1}
\ee
This system can be also expressed via the Fourier coefficients if we pass to the corresponding half-currents $x^i_{\pm}(u)$. Besides, we introduce the permutation relations
\be
x^a(u)\, R_{a\, i}^{b\,j}\, x_b(v)=\qq \, x_i(v)\, x^j(u)+\de_i^j \,\de(u-v),
\label{pe}
\ee
which are compatible with the  relations on the currents
$$x_i(u)=x_i^{+}(u)+x_i^{-}(u)\,\,{\rm and}\,\, x^j(u)=x^j_{+}(u)+x^i_{-}(u).$$

Thus, we have again constructed a quantum double $(A,B)$ of Fock type. The algebras $A$ and $B$ are generated by the Fourier coefficients of the currents  $x^k(u)$
and $x_k(u)$, which respectively play the role of the annihilation and creation operators. Using the definition $\de(u-v)=\sum_{p\in \Z}\frac{v^p}{u^{p+1}}$, we can express
the system (\ref{pe}) in terms of the Fourier coefficients:
$$
x^a[k] R_{a\, i}^{b\,j}\, x_b[l]=\qq \, x_i[l]\, x^j[k]+\de_i^j\,\de_{k+l}^1  \quad \Leftrightarrow \quad
x^a[k]\, x_b[l] =\qq\, x_i[l]\, x^j[k]\, \Psi_{j \,b}^{\,i\, a}+B_b^a\,\de_{k+l}^1.
$$

Now, define the counit $\varepsilon_A$  as above by setting
$$
\varepsilon_A(1_A)=1_{\C},\quad  \varepsilon_A(x^a[k])=0, \quad \forall\, a,\,k.
$$
This allows us to get the action of the algebra $A$ on $B$:
$$
x^a[k]\triangleright x_b[l]=B_b^a \de_{k+l}^1\,\,{\rm and\,\, so\,\, on}.
$$

We call the algebra, defined by the relations (\ref{chet}), (\ref{rela1}), and (\ref{pe}), the QD of Zamolodchikov-Faddeev type.

Now, introduce the matrix $L(u)=\|l_i^j(u)\|_{1\leq i,j, \leq N}$ where $l_i^j(u)=x_i(u)\, x^j(u)$.

\begin{proposition} The matrix $L(u)$ meets the following  relation
\rm
\be
R_{12}(u,v) L_1(u) R_{12} L_1(v)-L_1(v) R_{12} L_1(u)R_{12}(u,v)=(R_{12} L_1(u)-L_1(u) R_{12})\,\de(u-v).
\label{Yang}
\ee
\end{proposition}

\noindent
{\bf Proof.} Again, we use the matrix notation. Consider the first summand in the left hand side of (\ref{Yang}) and substitute the matrices $L(u)$ expressed  in terms of
currents $x_i(u)$ and $x^j(u)$:
\begin{eqnarray*}
&&R_{12}(u,v)x_{|1\rangle}(u)\underline{x^{\langle1|}(u) R_{12} x_{|1\rangle}(v)}x^{\langle 1|}(v) \stackrel{(\ref{pe})}{=}
q^{-1}\underline{R_{12}(u,v)x_{|1\rangle}(u)x_{|2\rangle}(v)}x^{\langle 2|}(u)x^{\langle 1|}(v)  \\
\rule{0pt}{5mm}&&+R_{12}(u,v)x_{|1\rangle}(u)x^{\langle 1|}(v)\,\delta(u-v)\stackrel{(\ref{chet})}{=} q^{-1}g(u,v)x_{|1\rangle}(u)x_{|2\rangle}(v)x^{\langle 2|}(u)x^{\langle
1|}(v)\\
\rule{0pt}{5mm}&& +R_{12}(u,v)x_{|1\rangle}(u)x^{\langle 1|}(v)\,\delta(u-v).
\end{eqnarray*}

In a similar way (the only difference is  using  (\ref{rela1}) instead of (\ref{chet})) we transform the second summand in the left hand side of (\ref{Yang}):
$$
L_1(v) R_{12} L_1(u)R_{12}(u,v) = q^{-1}g(u,v)x_{|1\rangle}(u)x_{|2\rangle}(v)x^{\langle 2|}(u)x^{\langle 1|}(v)+x_{|1\rangle}(v)x^{\langle 1|}(u)R_{12}(u,v)\,\delta(u-v).
$$

Taking the difference  of the above expressions we get the desired result:
\begin{eqnarray*}
&&\hspace*{-5mm}R_{12}(u,v) L_1(u) R_{12} L_1(v)-L_1(v) R_{12} L_1(u)R_{12}(u,v) = \\
\rule{0pt}{5mm}&&\hspace*{-5mm}(R_{12}(u,v)x_{|1\rangle}(u)x^{\langle 1|}(v) - x_{|1\rangle}(v)x^{\langle 1|}(u)R_{12}(u,v))\delta(u-v) =
(R_{12} L_1(u)-L_1(u) R_{12})\,\de(u-v).
\end{eqnarray*}
Note that the product $R(u,v)\delta(u-v)$ contains the ill-defined terms $\frac{(q-\qq) u}{u-v}\,\de(u-v)$, but they cancel each other
in the final expression.  As a consequence, the right hand side of (\ref{Yang}) contains only constant
matrix $R$.\hfill \rule{6.5pt}{6.5pt}

Concluding the paper we want to observe that the algebras analogous to those defined by (\ref{Yang}), but with vanishing right hand side, were introduced in \cite{GS1, GS2}
under the name of the braided Yangians or generalized Yangians of RE type. Similar algebras but with current braidings in the middle positions and depending of a charge
were considered in \cite{GS4}. Note that the algebras from \cite{GS4} are generalizations of those from \cite{RS}, corresponding to the $A_n$ series.

Nevertheless, the defining relations of the algebras from  \cite{GS4} cannot be expressed via the Fourier components. By contrast, the system (\ref{Yang}) can be written in terms
of these components. This is a consequence of the fact that we define  the permutation relations (\ref{pe}) in the QD above by means of a constant braiding (namely, a Hecke
symmetry), in contrast with the initial definition of the ZF algebras.

\end{document}